\documentclass[a4paper]{article}
\usepackage[margin=25mm]{geometry}
\usepackage{subfigure}
\usepackage{float}
\usepackage{graphicx}
\usepackage{amssymb}
\usepackage{mathrsfs}
\usepackage{fullpage}
\usepackage{amsmath}
\usepackage{graphicx,caption}
\usepackage{eucal}
\usepackage{lipsum}
\usepackage{amsfonts}
\DeclareMathAlphabet{\mathpzc}{OT1}{pzc}{m}{it}
\usepackage{epsfig}
\usepackage{adjustbox}
\allowdisplaybreaks
\usepackage[round]{natbib}
\usepackage{latexsym}
\usepackage[nodisplayskipstretch]{setspace}
\makeatletter
\def\ps@pprintTitle{%
\let\@oddhead\@empty
\let\@evenhead\@empty
\def\oddfoot\@empty
\let\evenfoot\@oddfoot}
\makeatother
\usepackage{verbatim}
\providecommand{\keywords}[1]
{
  \small	
  \textbf{\textit{Keywords---}} #1
}

\title{Numerical Optimization of  Loss System with Retrial Phenomenon in Cellular Networks} 

\author{ Vidyottama Jain$^{1}$, Raina Raj$^{1}$ and S. Dharmaraja$^{2}$  \\
        \small Central University of Rajasthan, Ajmer, India$^{1}$\\
       \small Indian Institute of Technology Delhi,  India$^{2}$\\
}

\date{}

\begin{document}
\maketitle

\begin{abstract}
In this study, we extend upon the model by Haring et al. [IEEE Trans. Veh. Technol. 50,  664-673 (2001)]  by  introducing retrial phenomenon in multi-server queueing system.   When at most $g$ number of guard channels are  available, it allows  new calls to join the retrial group. This retrial group is called orbit and can hold a maximum of $m$ retrial calls. 
The impact of retrial over certain performance measures is numerically investigated. The focus of this work is to construct  optimization problems  to determine the optimal number of channels, the optimal number of guard channels and the optimal orbit size. Further,  it has been emphasized that the proposed model with retrial phenomenon reduces the blocking probability of new calls in the system.
\end{abstract}


\section{Introduction}

\noindent Getting motivation from the  previously reported  research models, specifically,  \cite{gia} and  \cite{haring}, a multi-server queueing system with retrial phenomenon  is  studied here.  The retrial phenomenon  (\cite{dharma1}) is considered as when new calls are blocked due to non availability of idle channels and consequently join the orbit for retry. In this work,  such blocked new calls in the orbit will be referred as retrial calls. In a cellular network, termination of  new calls and handoff calls are crucial factors to determine the performance of the system. Therefore, dropping probability of handoff calls and blocking probability of new calls are the most  essential  performance measures for any cellular networks. In a realistic scenario, the service provider  will always  prefer  handoff calls over new calls. Henceforth, the service quality of handoff calls might be  improved by reserving   a group of guard channels (\cite{gurein})   for handoff calls.    Blocking probability of new calls has been a very important concept from 2G to today's 5G. Hence,  we address the problem of reducing  blocking probability of new calls by introducing retrial concept in cellular networks for a multi-server model with guard channels for handoff calls and retrial calls.

A lot of research has been carried out over the customer retrial phenomenon in a cellular network. Some of relevant studies are discussed here.
 \cite{marsan1} proposed a novel approximate multi-server model for the evaluation of call blocking probabilities in mobile cellular network taking into account retrial phenomenon. Another related work is \cite{gia}, in which   retrial model  with guard channel policy  was proposed, and a nearly recursive algorithm was derived  for providing  state probabilities.  \cite{sudhesh} studied transient analysis of a single-server retrial queueing system.  \cite{wang} analyzed a multi-server retrial queueing system without considering the guard channel policy.  They investigated the  behaviour of blocking probability and showed that   blocking probability is decreasing for retrial queues in case of multi-servers. Additionally, there are several  other techniques available in the literature   to evaluate the performance  of multi-server retrial queueing models (\cite{tien}, \cite{trivedi},  \cite{madan}, \cite{marson}, etc.). Recently  \cite{duc} proposed   a multi server retrial queueing model with two types of nonpersistent customers and studied their give up behaviour. They developed a numerically stable algorithm to compute the joint stationary distribution.

The most relevant work for our paper is \cite{haring},  which derived closed-form expressions for blocking probability of a new call and dropping probability of a handoff call.   \cite{haring} derived recursive formulae to compute  loss probabilities  for a multi-server queueing model and optimized  the  number of guard channels.  In this work,  an extension and generalization of \cite{haring} is proposed.   It  investigates the impact of a retrial on various performance measures including  blocking probability, dropping probability, mean number of busy channels and mean number of retrial calls. In this present study,  there is no closed form expression reported due to its complexity.  However, the  numerical results and optimization problems  presented here  consider \cite{haring}  as a particular case when there is no retrial.

This remainder of this paper is arranged as follows. Section 2 provides a brief  description of   the proposed  multi-server retrial queueing model and   elaborates this model mathematically.    Section 3 illustrates the numerical investigation of this multi-server retrial queueing model and  explains the impact of various parameters, i.e., orbit size, retrial rate,  number of total channels and number of guard channels, call arrival rate, etc., on the given performance measures. Further, Section 4 presents few optimization problems and their solution algorithms to determine the optimal values of the total number of channels, the  total number of guard channels and orbit size.   At last, discussion and future directions are provided in Section 5.

\section{Mathematical Model}

This work considers   a homogeneous cellular system  where each cell is served by a unique base station. Each base station consists of a finite number of total channels, say $c$, termed as channel pool.  This system can be modelled as a multi-server queueing model with retrial phenomena as depicted in Figure 1.

 Assume that the arrival pattern of the  new calls and the handoff calls  follows  the independent Poisson process with  the arrival rate $\lambda_{n}$ and $\lambda_{h}$ respectively.   Define the total arrival rate $\lambda= \lambda_{n}+\lambda_{h}$.
%Assume that the total arrival rate $\lambda= \lambda_{n}+\lambda_{h}$ follows  the independent Poisson process where $\lambda_{n}$ and $\lambda_{h}$, arrival rates  of new calls and handoff calls,  follows  the independent Poisson process.  }
 The  call duration  of new calls and handoff calls are i.i.d. random variables, each follows    independent exponential distribution with the rate  $\nu$. Note that   a handoff call is dropped when all channels are busy in the channel pool,  whereas  a new call joins the retrial group  when at least  $c-g$ channels are busy.  The  new call, which joins the retrial group, is referred as retrial call in this work.    Here, it is assumed that this retrial group, called orbit,  has a capacity of finite size, say $m$.  Once the orbit is occupied by  $m$  retrial calls, the new call will not be able to enter into the system and therefore, it will be blocked. In this model, a very well known guard channel policy is considered.   Under this policy,  $g$ number of channels are reserved for handoff calls and retrial calls.  In the proposed model, it is considered  that a retrial call   retries at random intervals and in random order.  It obtains the service  with  probability $p$     or leaves the system forever  without obtaining the service with probability $1-p$.  The inter-retrial time between the  retrial calls  is exponentially distributed with rate $\mu_r$.

  %On retrial, with retrial probability $p$, a retrial call from  the orbit gets connected  immediately.    In this model, 
\begin{figure}
\centering
\caption{A multi-server queueing model with retrial phenomenon.}
\includegraphics[width=4.3in, height=3.8in]{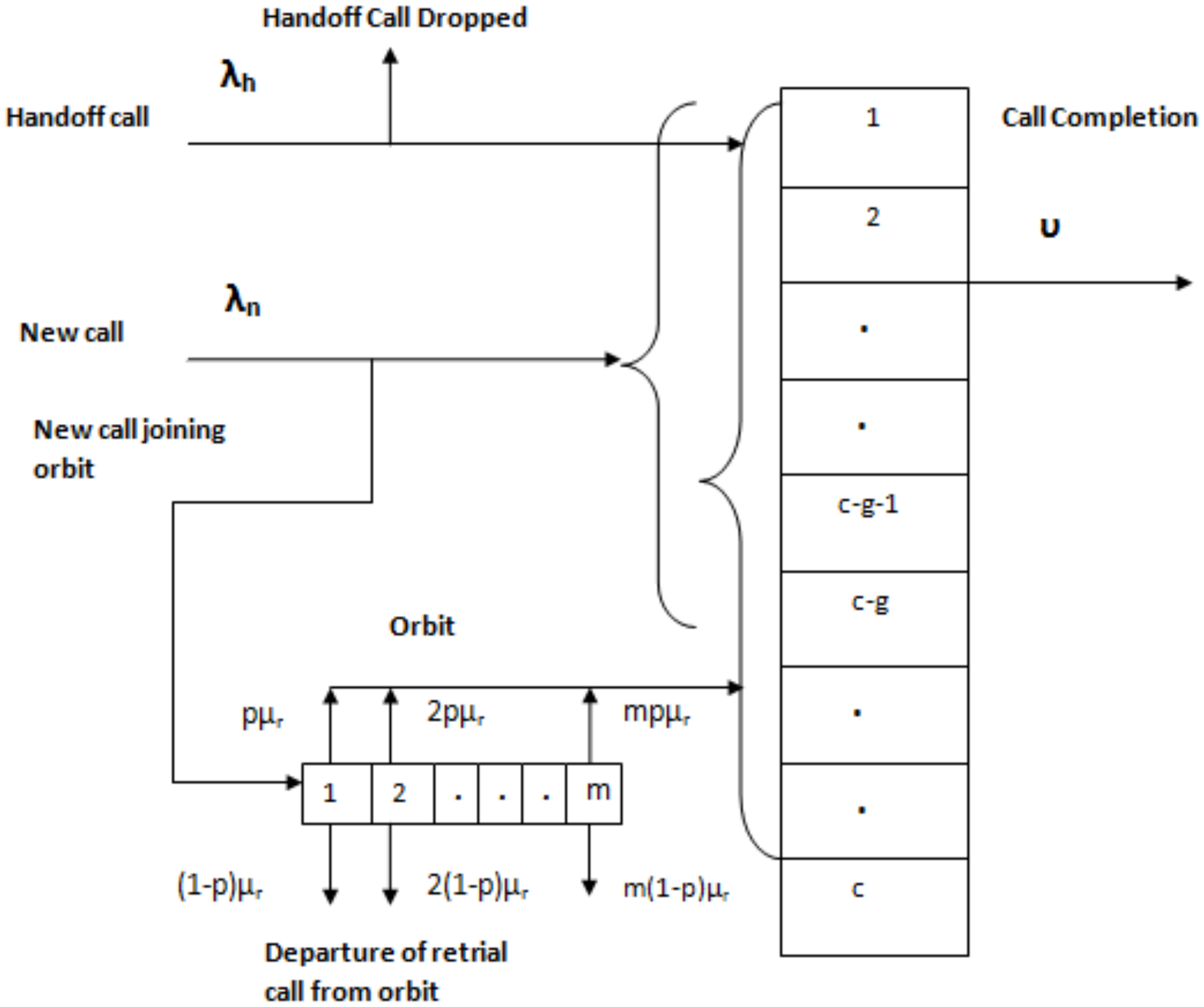} 
\label{fig:fig_sim1}
\end{figure} 

\begin{figure}
\centering
\caption{State transition diagram}
\includegraphics[width=5in]{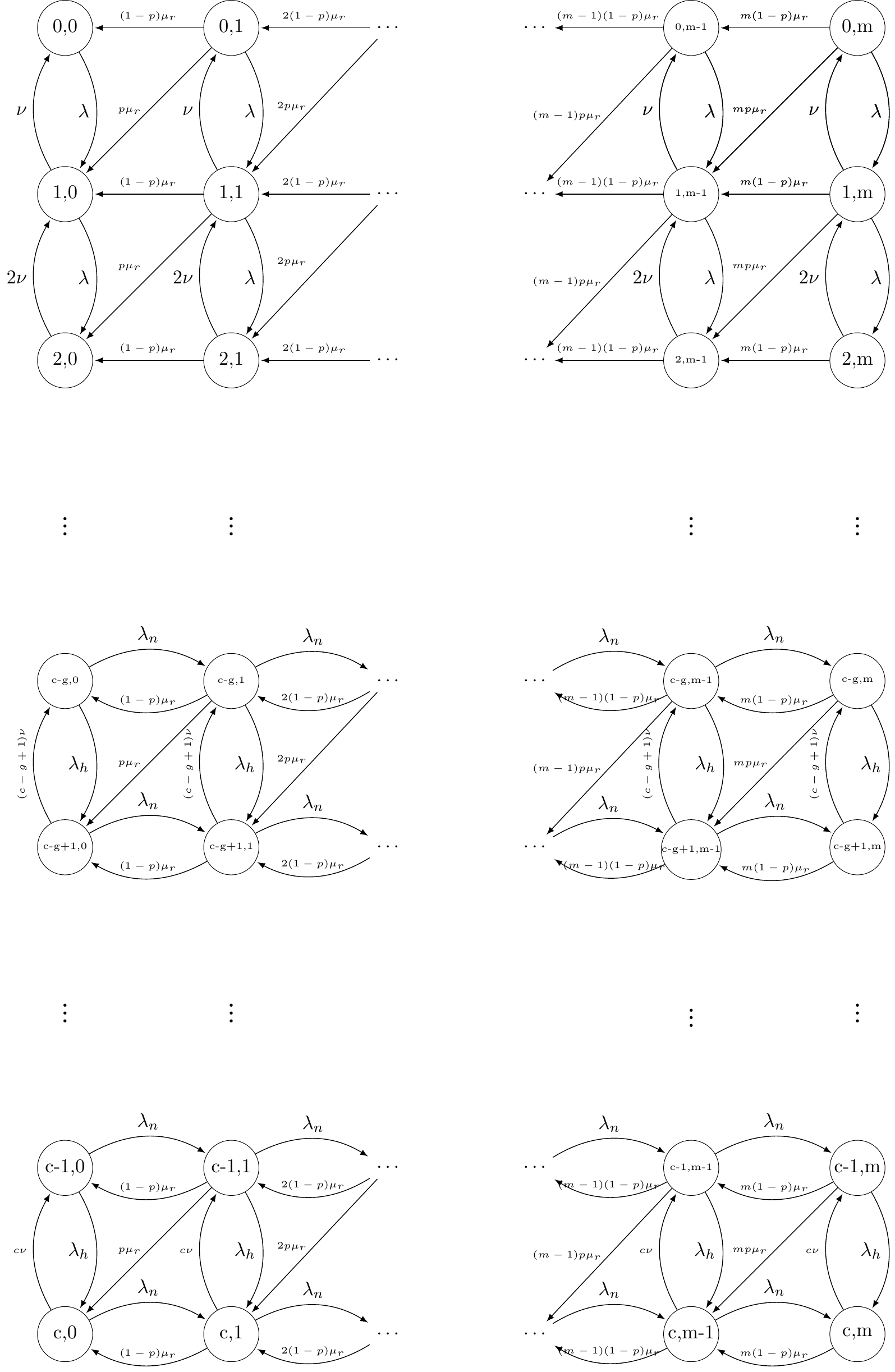}
\label{fig:fig_2}
\end{figure}

 Let us analyze the underlying  stochastic process for the proposed  queueing model.  Suppose $\Xi(t)$ represents the number of busy channels at time $t$ and $\Theta(t)$ defines the number of blocked new calls in the orbit  at time $t$. Based on the model description,   the stochastic process $\{ (\Xi(t), \Theta(t)): t \geq 0 \}$ can be modelled as a quasi-birth-and-death (QBD) process  with a state space  $S = \{(j,k) : 0 \leq j \leq c ;~0 \leq k\leq m  \}$. 
Here $j$ denotes the number of busy channels at time $t$ and $k$ denotes the number of retrial calls in the orbit  at time $t$. Figure \ref{fig:fig_2} presents the state transition diagram for the proposed model.

Since the underlying stochastic process is ergodic, the stationary distribution of the system exists and is independent of initial distribution.  Let  $P_{j,k}$ be defined as 
\begin{align}
\nonumber P_{j,k} &= \displaystyle  \lim_{t \rightarrow \infty} P(\Xi(t)=j, \Theta(t)=k); ~0\le j \le c,~0\le k\le m. 
\end{align}

Then,  the steady state equations  can be written as follows:
\begin{align}
& \lambda  P_{0,0}  =    \nu P_{1,0} + (1-p)\mu_{r} P_{0,1},\\
&(\lambda + k\mu_{r} ) P_{0,k}  =    \nu P_{1,k} + (k+1)(1-p)\mu_{r} P_{0,k+1};  \quad 1 \leq k \leq m-1,\\
&(\lambda + m\mu_{r} ) P_{0,m}  = \nu P_{1,m} , \\
&\nonumber (\lambda + k\mu_{r} + j\nu) P_{j,k}  =   \lambda P_{j-1,k}  + (j+1)\nu P_{j+1,k} +  (k+1)p \mu_{r} P_{j-1, k+1}    \\ 
& \nonumber\hspace{3.3cm}  +   (k+1)(1-p)\mu_{r} P_{j,k+1}; \\& \hspace{3.3cm} 1 \leq j \leq c-g-1,   \quad 0 \leq k \leq m-1, \\
& (\lambda + m\mu_{r} + j\nu) P_{j,m}  =   \lambda P_{j-1,m}  + (j+1)\nu P_{j+1,m}; \quad 1 \leq j \leq c-g-1,  \\
& \nonumber(\lambda +  (c-g)\nu) P_{c-g,0}   =  \lambda P_{c-g-1,0}  + (c-g+1)\nu P_{c-g+1,0} \\ & \hspace{3.2cm} + p \mu_{r} P_{c-g-1,1}  + (1-p)\mu_{r} P_{c-g,1}, \\
&\nonumber(\lambda + k\mu_{r} + (c-g)\nu) P_{c-g,k}   =    (c-g+1)\nu P_{c-g+1,k} +  \\ & \hspace{4.2cm} \nonumber \lambda P_{c-g-1,k} +  (k+1)p\mu_{r} P_{c-g-1,k+1}   \\  & \hspace{4.1cm} +  (k+1)(1-p)\mu_{r} P_{c-g,k+1};   \quad 1 \leq k \leq m-1,\\
&(\lambda_h + m\mu_{r} + (c-g)\nu) P_{c-g,m} =   \lambda P_{c-g-1,m}  + (c-g+1)\nu P_{c-g+1,m},\\
&\nonumber(\lambda  + j\nu) P_{j,0}  =  \lambda_{h} P_{j-1,0}   + (j+1)\nu P_{j+1,0}  + p \mu_{r} P_{j-1,1} \\ & \hspace{2.2cm}   +  (1-p)\mu_{r} P_{j,1};    c-g+1 \leq j \leq c-1,\\
&\nonumber(\lambda + k\mu_{r} + j\nu) P_{j,k}  =  \lambda_{h} P_{j-1,k}  + \lambda_{n} P_{j-1,k-1} + (j+1)\nu P_{j+1,k} \\
   & \hspace{3.2cm} \nonumber + (k+1)p \mu_{r} P_{j-1,k+1}    +   (k+1)(1-p)\mu_{r} P_{j,k+1}; \\ & \hspace{3.2cm}   c-g+1 \leq j \leq c-1,  \quad 1 \leq k \leq m-1,\\
&\nonumber(\lambda_h + m\mu_{r} + j\nu) P_{j,m}  =  \lambda_{h} P_{j-1,m}  + \lambda_{n} P_{j-1,m-1} + (j+1)\nu P_{j+1,m} ; \\ &  \hspace{3.2cm}  c-g+1 \leq j \leq c-1, \\
&(\lambda_{n} + c\nu)P_{c,0}  =  \lambda_{h} P_{c-1,0}  + p \mu_{r} P_{c-1, 1} + (1-p)\mu_{r} P_{c,1},\\
&\nonumber(\lambda_{n} + k(1-p)\mu_{r} + c\nu)P_{c,k}  =  \lambda_{h} P_{c-1,k}+ \lambda_{n} P_{c,k-1}  + (k+1)p \mu_{r} P_{c-1, k+1} \\ & \hspace{4cm}  +  (k+1)(1-p)\mu_{r} P_{c,k+1}; \quad 1 \leq k \leq m-1,\\
&( m(1-p)\mu_{r} + c\nu)P_{c,m} =  \lambda_{h} P_{c-1,m}+ \lambda_{n} P_{c,m-1}. 
\end{align}

Further, with these steady state equations, the infinitesimal generator matrix $Q$ for the QBD $\{ (\Xi(t), \Theta(t)): t\geq 0\}$ is constructed as follows:
\[
Q =
\begin{pmatrix}
  Q_{0,0} & Q_{0,1} & 0 &  &  &  \\
      Q_{1,0}& Q_{1,1}& Q_{1,2} &  &  &  \\
      0  & Q_{2,1} & Q_{2,2} &  &  &   \\
      0 & 0 & Q_{3,2}  &  &  &               \\
      0 & 0 & 0  &  &  &               \\
       &  & \ddots & \ddots & \ddots & \\
       &  &    &  &         & Q_{c-1,c}\\
        &  &    &  &      Q_{c,c-1}   & Q_{c,c}
                                  
\end{pmatrix}.\]

Each element of matrix $Q$ is a block matrix of order $(m+1) \times (m+1)$. More specifically, this tridiagonal matrix $Q$  is described as:
\begin{itemize}
\item The elements of the upper diagonal, denoted as block matrix $Q_{l,l+1}$ where $0\leq l \leq c-1$, shows the transition due to the arrival of a (handoff or new) call or joining the orbit by a new  call for retrial. Therefore, structure of $Q_{l,l+1}$ is as follows:

\[
Q_{l,l+1} = 
\begin{pmatrix}
\Lambda  & 0 & 0 &  & \\
p\mu_{r} & \Lambda & 0  & &  \\
0 & 2p\mu_{r}  & \Lambda & &\\
0 & 0 & 3p\mu_{r}  &  & &\\
 &   & \ddots & \ddots &\\
 &  &  &mp\mu_{r} & \Lambda

\end{pmatrix}\]

where $\Lambda =  \begin{cases}
     \lambda=\lambda_h+\lambda_n, ~~~~~~~ 0\leq l \leq c-g-1,\\    
     \lambda_h,~~~~~~~~~~~~~~~~~~~~c-g\leq l \leq c-1.   
\end{cases} $
                   
\item The elements of the main diagonal, denoted as $Q_{l,l}$  where $0\leq l \leq c$, represents the transitions either due to the departure of a call after its completion or arrival of a new call or departure of a retrial call without getting the connection. Consequently, the block matrix $Q_{l,l}$ is constructed as 

\[Q_{l,l} = 
\begin{pmatrix}
b_0& \Lambda' & 0 & 0 &  & \\
d_1 & b_1 & \Lambda' & 0 & &  \\
0 &d_2 & b_2 & \Lambda' &  & \\
 &  & \ddots & \ddots & \ddots & \\
 &  &  & \ddots & \ddots & \Lambda' \\
&  &  &  & d_{m} & b_{m}

\end{pmatrix} 
\]
where $\Lambda' =  \begin{cases}
     0, ~~~~~~~ 0\leq l \leq c-g-1,\\    
     \lambda_n,~~~~~~c-g\leq l \leq c.  
\end{cases} $

$\quad$ $\quad$ $b_k = -(\lambda + l\nu + k\mu_{r})$; $0\leq k \leq m$,
 
and 
  
$\quad$ $\quad$ $d_k =  k(1-p)\mu_{r}$;  $1\leq k \leq m$.

\item The elements of the lower diagonal, denoted as $Q_{l,l-1}$ where $1\leq l \leq c$,  exhibits the transitions due to the completion of a call (new or handoff). Hence, these block matrices $Q_{l,l-1}$ are presented as 

\[Q_{l,l-1} = 
\begin{pmatrix}
 l\nu& 0& &   \\
0 &l\nu &  &  \\
 &  & \ddots   &   \\
&  &   &  l\nu

\end{pmatrix}.\]

\end{itemize}

Finally, the steady state equations (1)-(14)   can be expressed in matrix form as $\Pi Q = 0$,  where
 \[\Pi = (\Pi_{0}, \Pi_{1}, \Pi_{2}, \ldots, \Pi_{c})\]
 and 
 \[\Pi_{j} = (P_{j,0}, P_{j,1},\ldots,P_{j,m}), ~~ 0 \leq j \leq c,\]
 
\noindent  with the normalization condition $\Pi e = 1$ where $e$ is a unit vector. Here  $Q$ is a complicated and highly structured matrix, therefore,  it is difficult to obtain  a compact  analytical form for $\Pi$.  This system is solved by applying the  direct method to compute   steady state probabilities.

\subsection{Performance Measures}

With the steady state probabilities,  the relevant  performance measures for  the proposed model  are  represented in this section.
\begin{itemize}
\item The  probability of new calls being blocked  is computed by:
   \[ P_b =  \sum_{j=c-g}^{c} P_{j,m}.\]
\item The probability of handoff call being dropped can be calculated as:
 \[P_d = \sum_{k=0}^{m} P_{c,k}.\]
%\item The probability of first attempt fresh new calls being blocked is given by:
%\[ P_{bf} =  \sum_{j=c-g}^{c} P_{j,0}.\]
\item The mean number of busy channels in the channel pool can be evaluated  by:
\[ M_b = \sum_{j=1}^{c}\sum_{k=0}^{m} j P_{j,k}.\]
\item The mean number of repeated retrial calls in the orbit can be obtained by:
\[ M_o = \sum_{j=1}^{c}\sum_{k=0}^{m} k P_{j,k}.\]
\item The mean number of new calls and handoff calls in the system can be evaluated by:
\[M_s = M_b + M_o.\]
\end{itemize}

\section{Numerical Illustration}

The goal of this section is to analyse numerical results received after implementing  the mathematical model presented in Section 3.  This analysis has two main objectives: first objective is to present that this work  is an extension of the work provided by  \cite{haring}, and second objective is to examine how the  retrial phenomenon affects  the system performance. Therefore, for illustration purpose, we set the parameters  as  $\lambda=\lambda_{h} +\lambda_{n}= 80$, $\lambda_{n} = 40$, $\nu = 1$, $p = 0.8$  and $\mu_r = 0.5$.    Here, values of all these parameters are same as mentioned in   \cite{haring}  except retrial probability $p$ and retrial rate $\mu_r$. It is important to note that though  the number of channels $c$, number of guard channels $g$ and orbit size $m$ are  integers, yet these parameters are considered here as real numbers in order to plot figures and to analyze the results.

\begin{figure*}
\centering
\caption{Dependence of loss probabilities over  $c$ for $m=0$.}
\subfigure[Dependence of $P_b$ over $c$ for $m=0$. ]
 {\includegraphics[trim= 1.5cm 0.1cm 3cm 0.5cm, clip=true, height = 6cm,width = 0.54\textwidth]{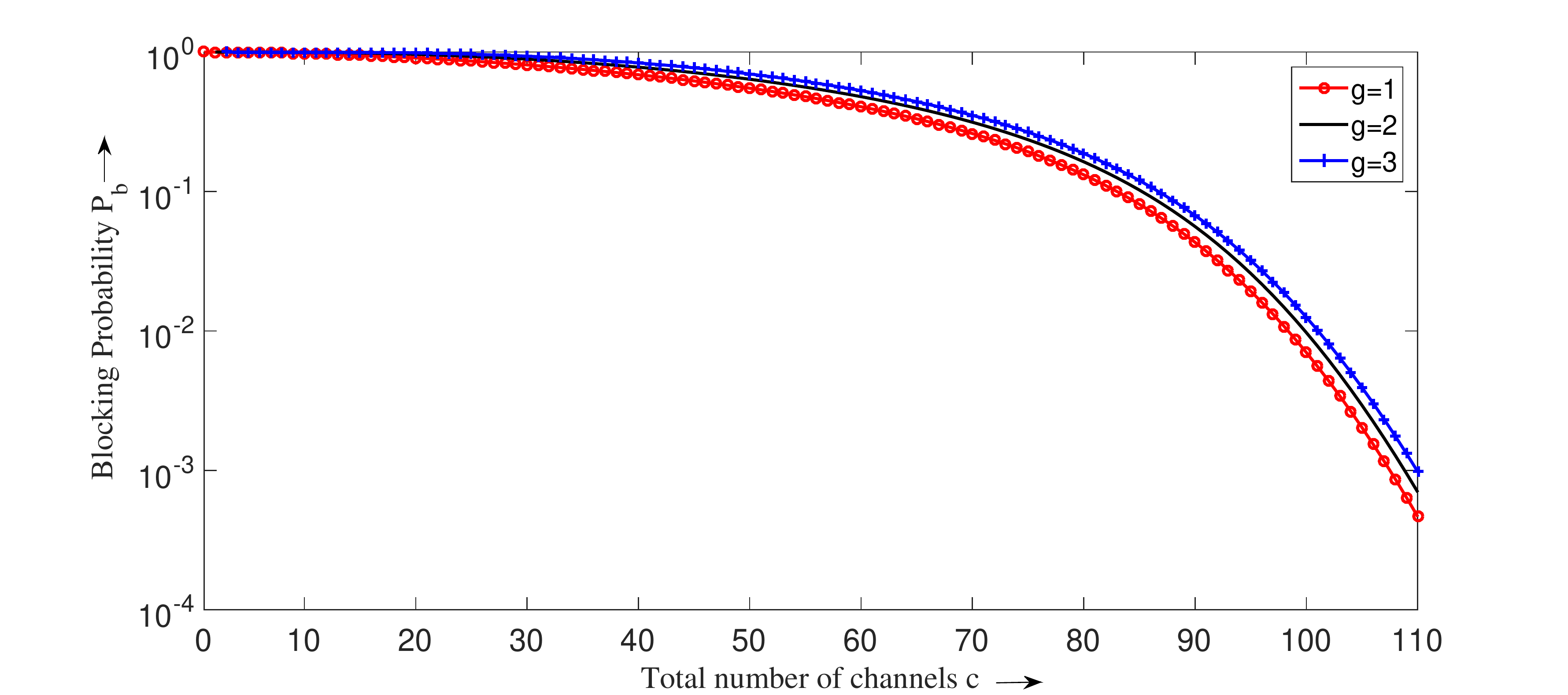} 
  \label{fig:Pb_N}
  }%
\subfigure[Dependence of $P_d$ over $c$ for $m=0$.]
{\includegraphics[trim= 2cm 0.1cm 1cm 0.3cm, clip=true, height = 6cm,width = 0.56\textwidth]{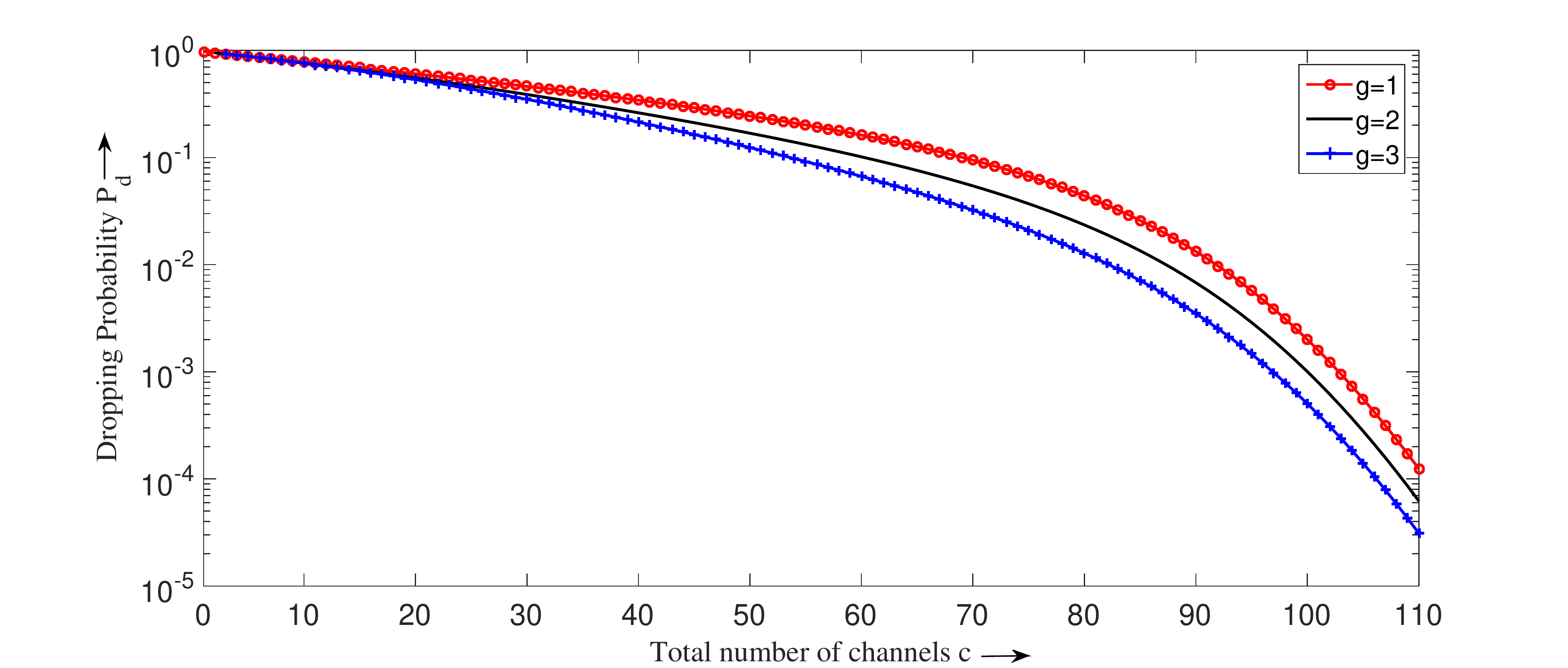} 
\label{fig:Pd_N}
}%

\label{fig:PbPd_N}
\end{figure*}

\begin{figure*}
\centering
\caption{Dependence of $P_b$ and $P_d$ over $c$ and $g$ for $m=0$.}
{\includegraphics[trim= 1.5cm 0.1cm 1cm 0.1cm, clip=true, height = 6.5cm,width = 1\textwidth]{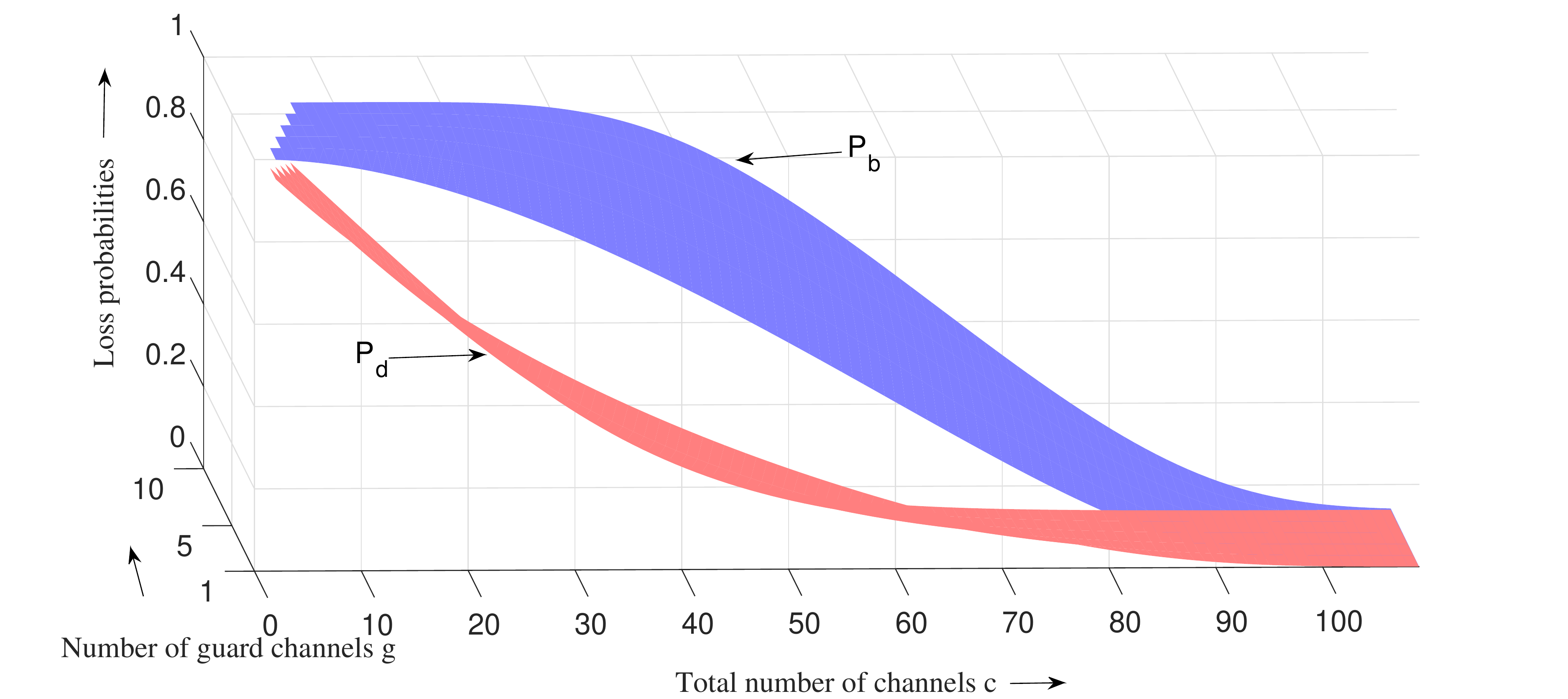}

\label{fig:constrained_graphs}
}%
\end{figure*}

After performing the numerical illustration for  $g=1,2,3$ and  $m=0$, Figure \ref{fig:PbPd_N}\subref{fig:Pb_N}  and Figure \ref{fig:PbPd_N}\subref{fig:Pd_N}    plot  loss probabilities (the blocking probability  $P_b$ and the dropping probability  $P_d$)  with respect to  $c$. We observe  by these figures  that  for no retrial, i.e., for $m=0$,  $P_b$ and  $P_d$ decreases with $c$ for a fixed value of $g$,    $P_b$ increases with $g$ for a fixed value of $c$, and $P_d$ decreases with $g$ for a fixed value of $c$. In this regard, Figure   \ref{fig:constrained_graphs} clearly displays the dependence of loss probabilities  over $c$ and $g$ simultaneously for  $m=0$.
 These properties of $P_b$ and $P_d$ match well with the  properties provided by  \cite{haring}. Therefore, we declare that, in case of no retrial,  this work provides same properties of $P_b$ and $P_d$ as reported by  \cite{haring}.

\begin{figure*}
\centering
\caption{Dependence of loss probabilities over  $g$ for $c=100$. }
\subfigure[Dependence of $P_b$ over $g$ for $c=100$.]
{\includegraphics[trim= 2.1cm 0.01cm 3.1cm 0.001cm, clip=true, height = 6cm,width = 0.54\textwidth]{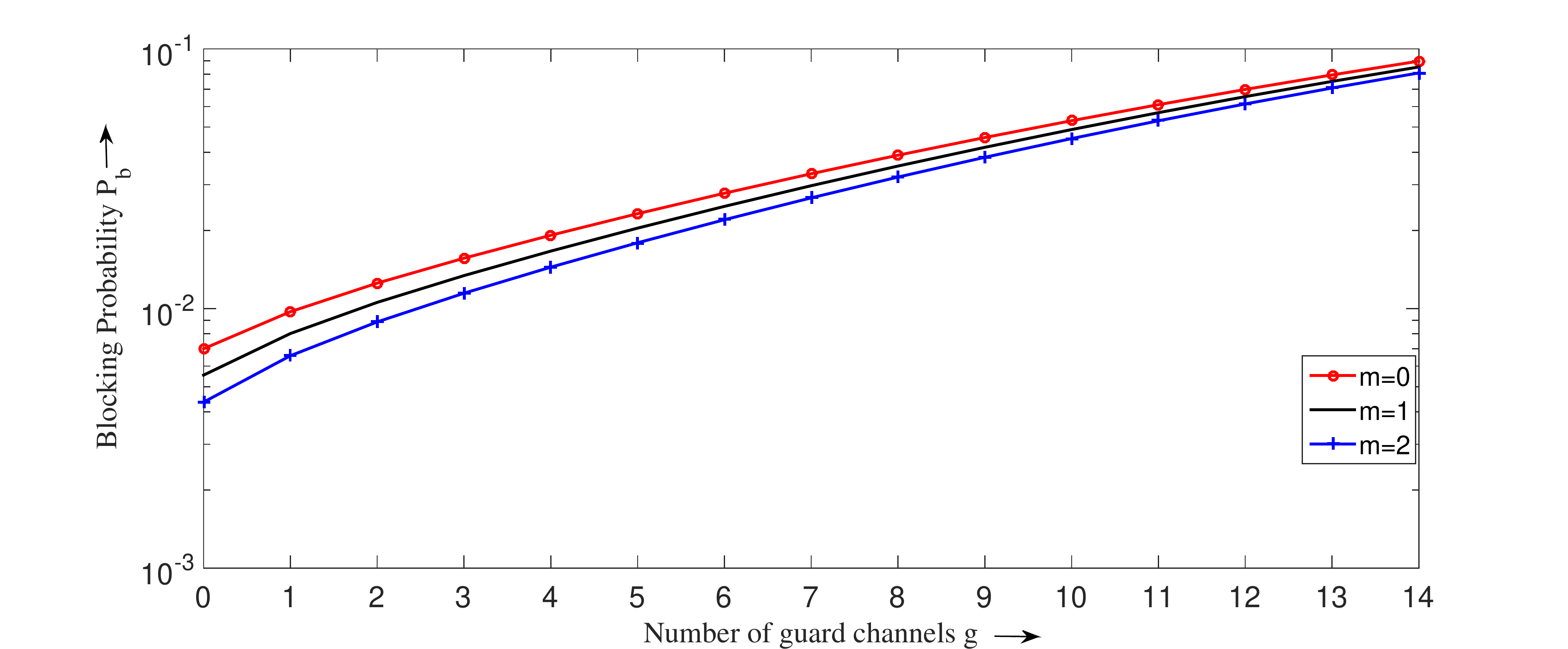}  
 \label{fig:Pb_g}
  }%
\subfigure[Dependence of $P_d$ over $g$ for $c=100$.]
{\includegraphics[trim= 2.1cm 0.01cm 0.8cm 0.001cm, clip=true, height = 6cm,width = 0.58\textwidth]{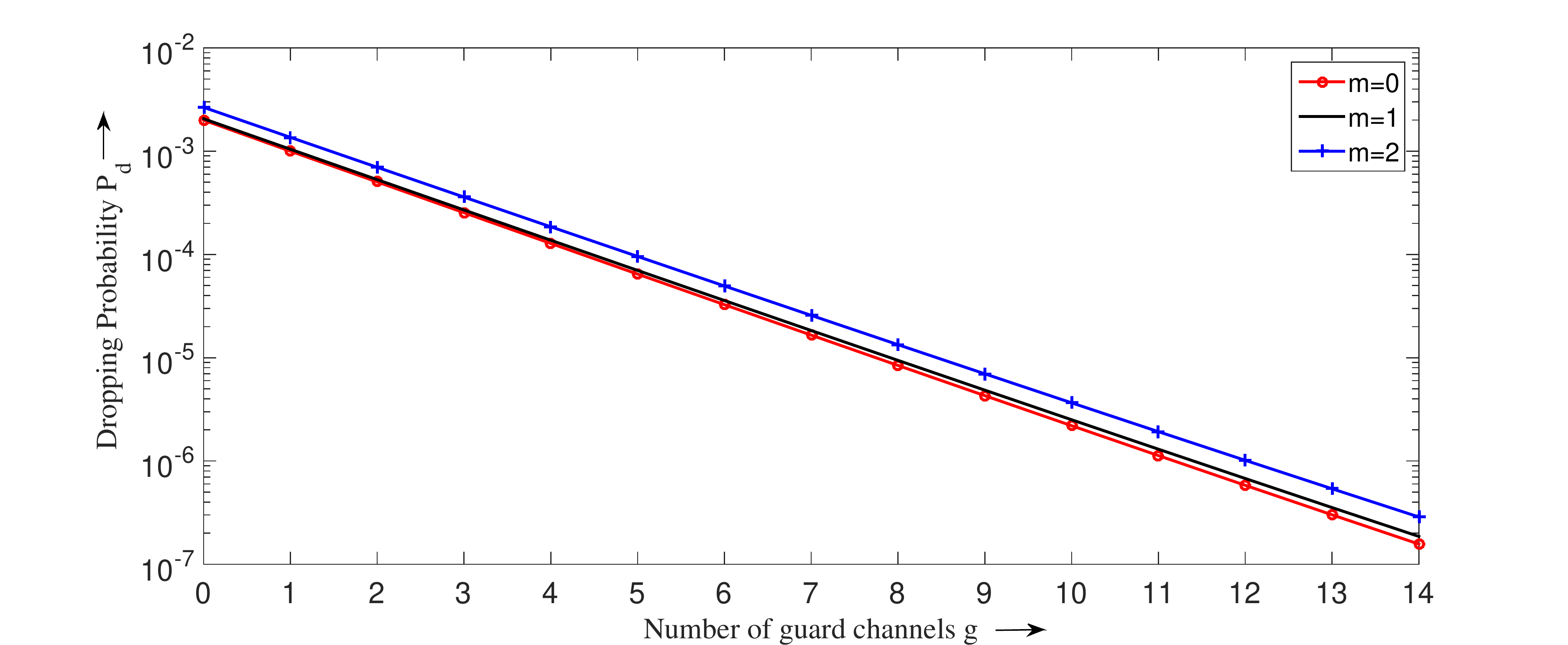} 
\label{fig:Pd_g}
  }%

\label{fig:PbPd_g}
\end{figure*}

The  remaining part of this section demonstrates the effect of  retrial phenomenon over the various  performance measures given in previous section. For $m>0$, Figure \ref{fig:PbPd_g}\subref{fig:Pb_g} and Figure \ref{fig:PbPd_g}\subref{fig:Pd_g} exhibit  that $P_b(g)$ is an increasing function and  $P_d(g)$ is a decreasing function for a fixed value of $c$.  Figure \ref{fig:3D_Pb(N,g)} \subref{fig:3D_Pb(N,g)_1} and Figure \ref{fig:3D_Pb(N,g)} \subref{fig:3D_Pd(N,g)_2} explore   that by increasing the orbit size, i.e., $m$,  $P_b$ decreases and $P_d$ increases.  Similar results can be obtained by examining  Figure \ref{fig:PbPd_m}\subref{fig:Pb_m} and Figure \ref{fig:PbPd_m}\subref{fig:Pd_m}. 

%===================================    

\begin{figure*}
\centering
\caption{Dependence of loss probabilities over $c$ and $g$ for $m=1, 4, 10, 20$. }
\subfigure[Dependence of $P_b$ over $c$ and $g$ for $m=1, 4, 10, 20$.]
{\includegraphics[trim= 1.5cm 0.1cm 1.8cm 0.8cm, clip=true, height = 6cm,width = 0.54\textwidth]{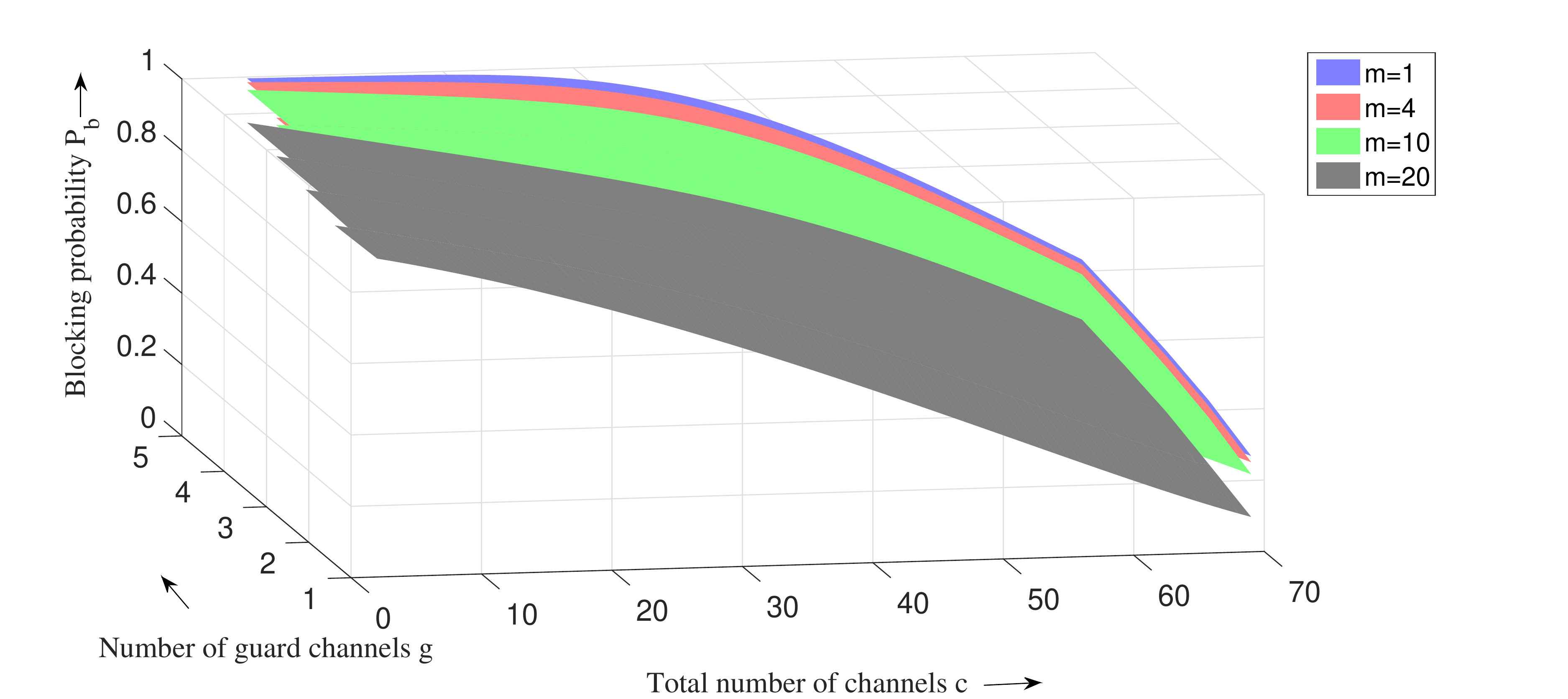} 
 \label{fig:3D_Pb(N,g)_1}
  }%
\subfigure[Dependence of $P_d$ over $c$ and $g$ for $m=1, 4, 10, 20$.]
{\includegraphics[trim= 1.1cm 0.1cm 1.8cm 0.6cm, clip=true, height = 6cm,width = 0.58\textwidth]{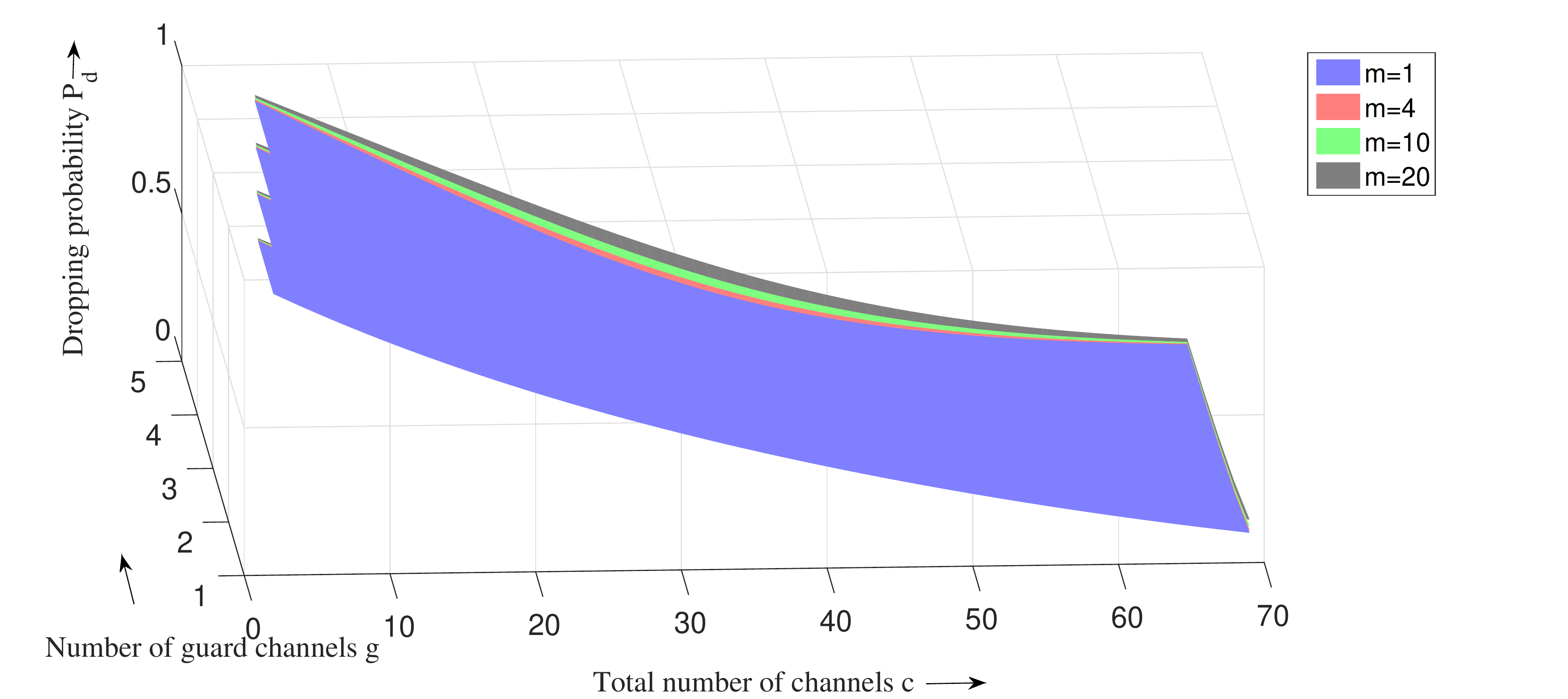}
 \label{fig:3D_Pd(N,g)_2}
  }%

\label{fig:3D_Pb(N,g)}
\end{figure*}

%======================================

\begin{figure*}
\centering
\caption{Dependence of loss probabilities over $m$ for $c=100$.}
\subfigure[Dependence of $P_b$ over $m$ for $c=100$.]
{
\includegraphics[trim= 1.5cm 0.1cm 3cm 0.8cm, clip=true, height = 6cm,width = 0.5\textwidth]{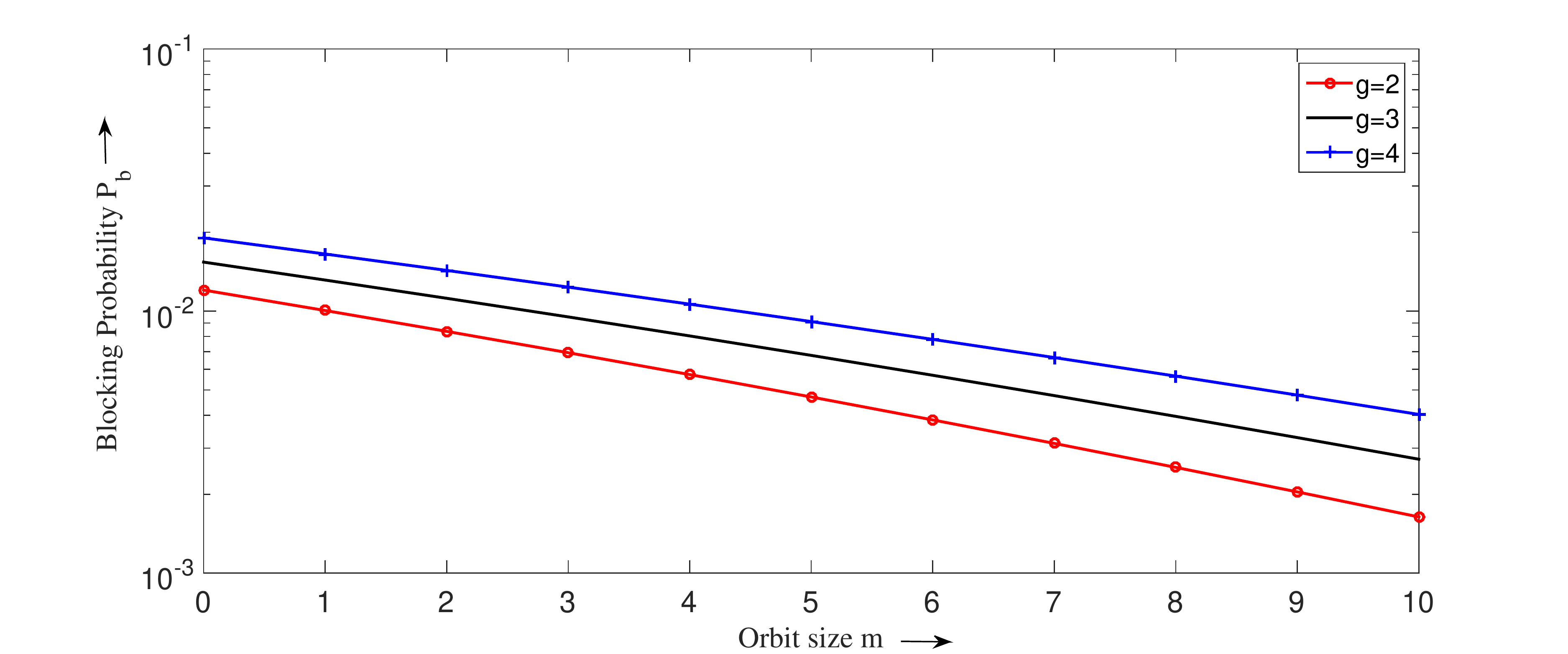}  
\label{fig:Pb_m}
}%
\subfigure[Dependence of $P_d$ over $m$ for $c=100$.]
{\includegraphics[trim= 2cm 0.2cm 1cm 0.5cm, clip=true, height = 6cm,width = 0.5\textwidth]{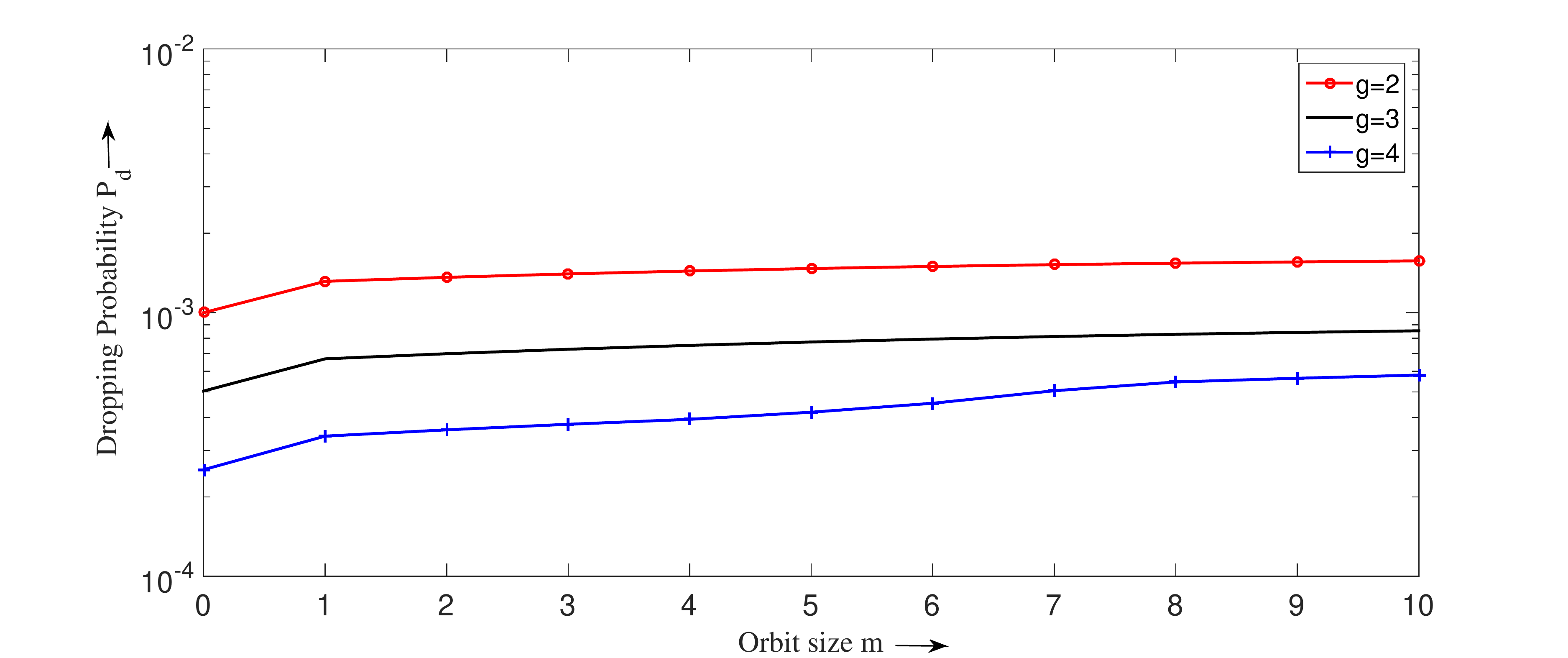} 
\label{fig:Pd_m}
}%

\label{fig:PbPd_m}
\end{figure*}

\begin{figure*}
\centering
\caption{Loss probabilities as a function of retrial rate $\mu_r$ for $c=100$, $m=5$ and $g=5$.}
{\includegraphics[trim= 1.5cm 0.2cm 3cm 0.1cm, clip=true, height = 6cm,width = 0.5\textwidth]{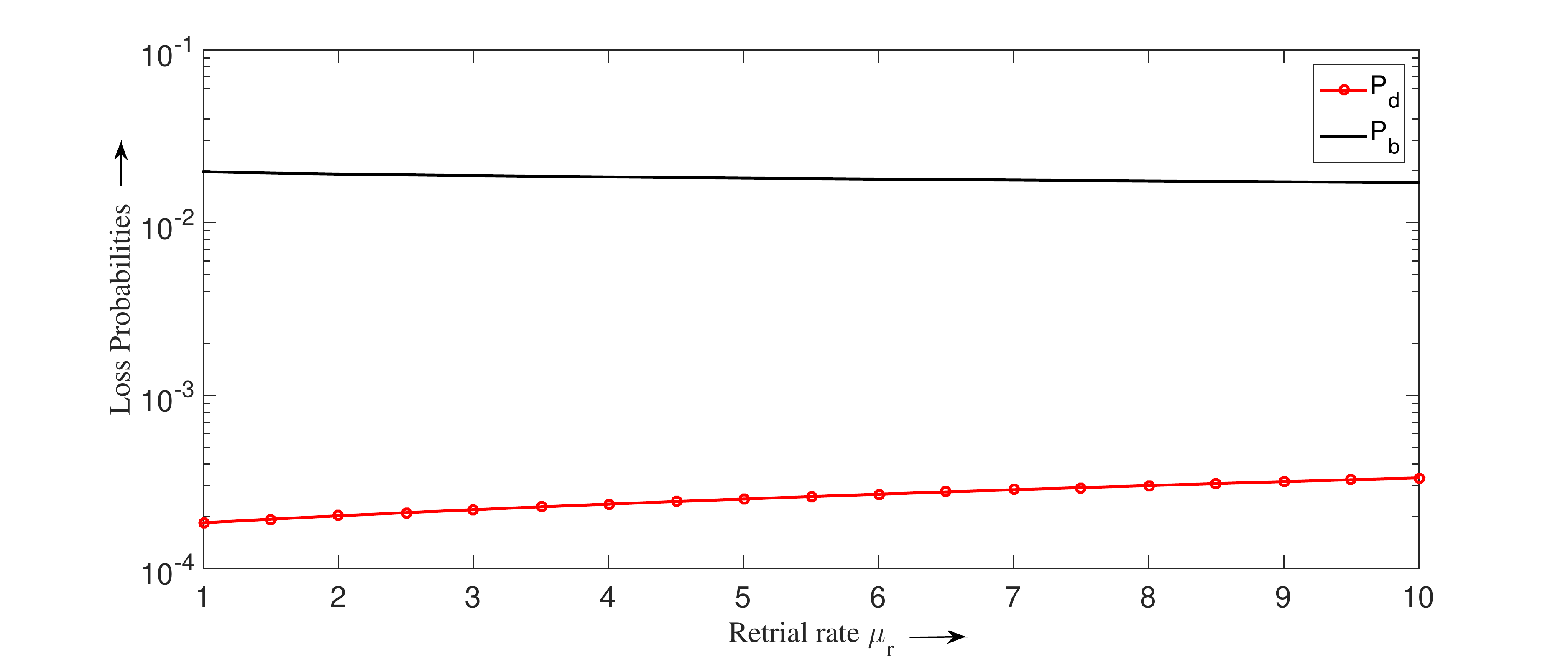}

\label{fig:loss_prob_retrial_rate}
}%
\end{figure*}

Figure \ref{fig:loss_prob_retrial_rate} reflects the impact of retrial rate on the loss probabilities. It is obvious that with the increased  retrial rate, the blocking probability is gradually reduced as more retrial calls get connected. With the similar reasoning, increment in retrial rate will increase the dropping probability $P_d$.

\begin{figure*}
\centering
\caption{Dependence of loss probabilities over  call arrival rates for $c=100$ and $m=5$.}
\subfigure[Dependence of $P_b$ over $\lambda_h$ for $c=100$ and $m=5$.]
 {\includegraphics[trim= 1.5cm 0.2cm 3cm 0.2cm, clip=true, height = 6cm,width = 0.5\textwidth]{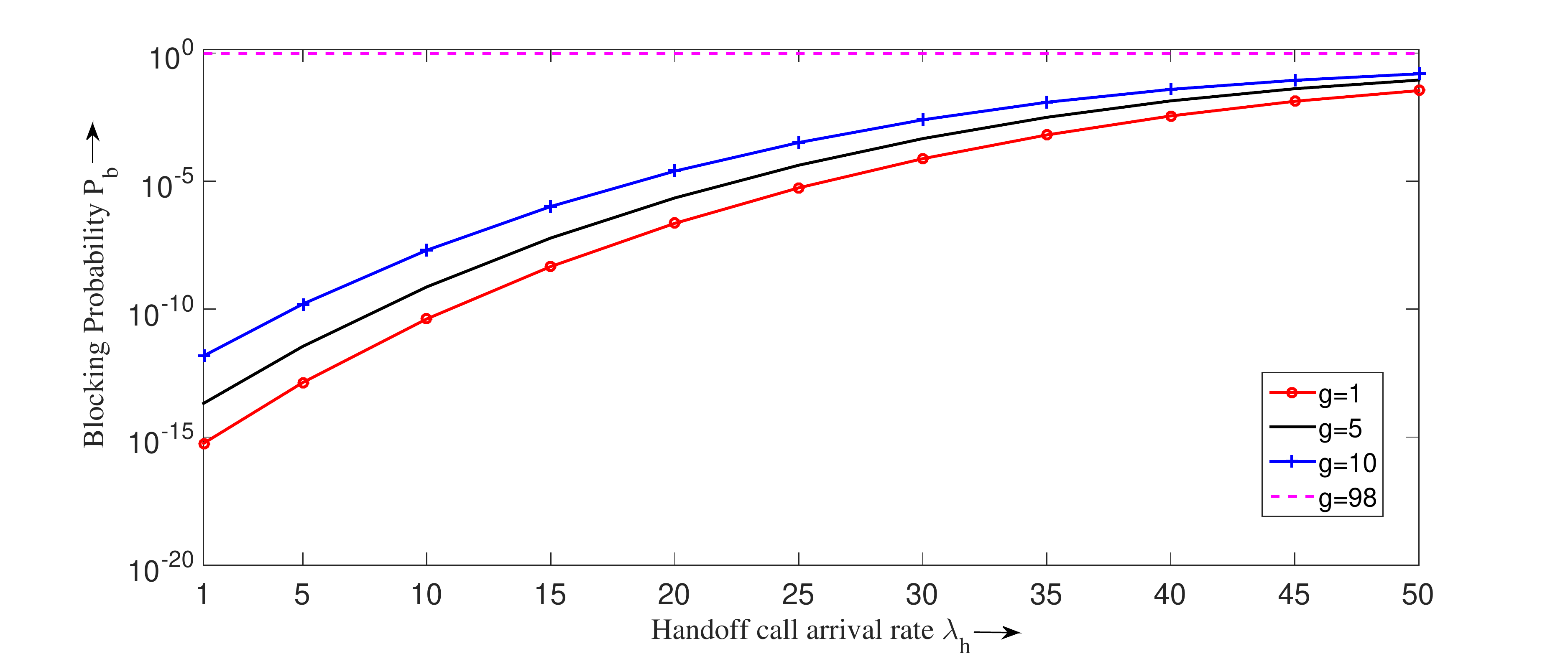}  
\label{fig:Pb_lambda}
}%
\subfigure[Dependence of $P_d$ over $\lambda_n$ for $c=100$ and $m=5$.]
 {\includegraphics[trim= 2cm 0.2cm 1cm 0.2cm, clip=true, height = 6cm,width = 0.5\textwidth ]{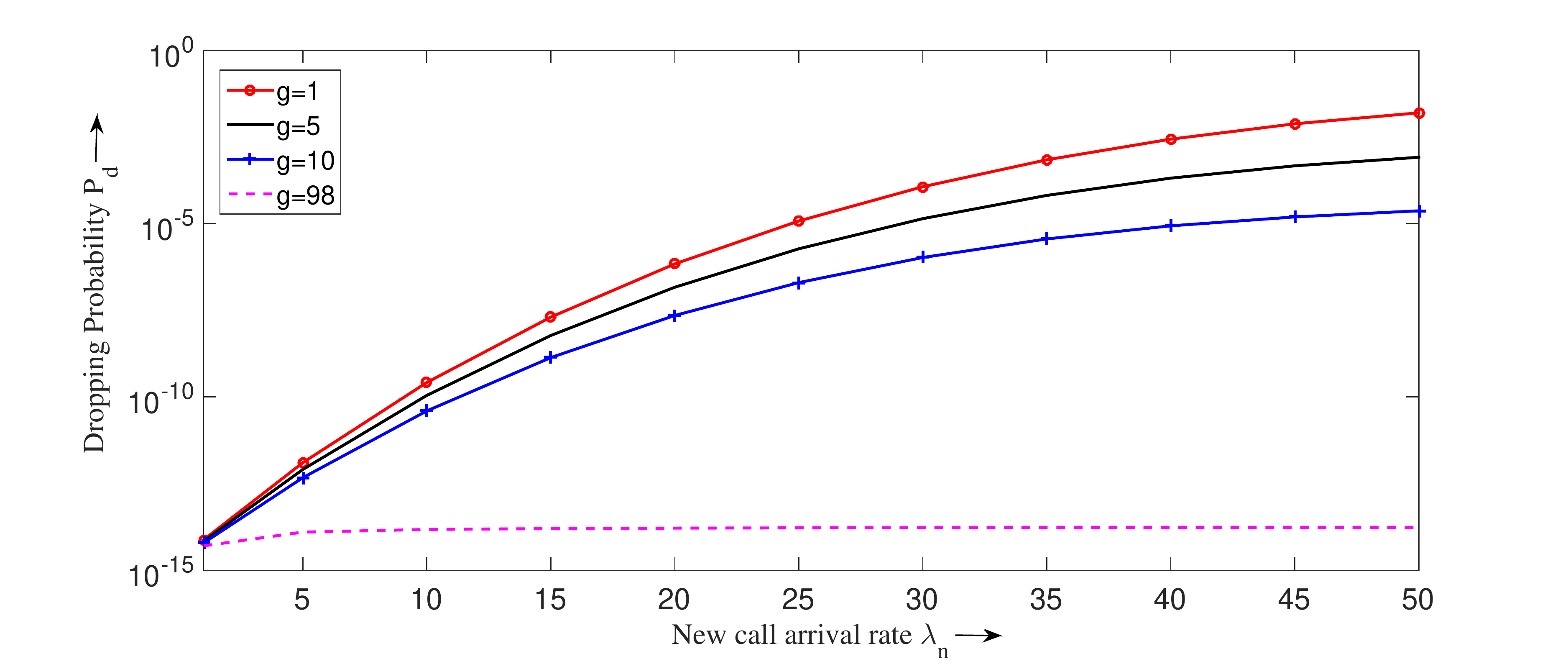} 
\label{fig:Pd_lambda}
}%

\label{fig:PbPd_lambda}
\end{figure*}

%##############################
Additionally, in case of retrial phenomenon,  we  can intuitively  perceive  that the increment of handoff call arrival rate will increase  $P_b$ and  the increment of  new call arrival rate will  increase $P_d$ for a fixed value of $g$. It can also be realized that for larger values of $g$, i.e. almost equal to $c$,  the arrival rate of new calls should have a negligible impact over  $P_d$.  Such effect of  $\lambda_h$ on  $P_b$ and   $\lambda_n$ on  $P_d$ are displayed by  Figure \ref{fig:PbPd_lambda}\subref{fig:Pb_lambda}  and Figure \ref{fig:PbPd_lambda}\subref{fig:Pd_lambda}. 

\begin{figure*}
\centering
\caption{Dependence of $M_b$ and  $M_o$ over  $c$ for $m=5$ and $g=5$.}
\subfigure[Dependence of $M_b$ over   $c$ for $m=5$ and $g=5$.]
{\includegraphics[trim= 1.5cm 0.15cm 3cm 0.5cm, clip=true, height = 6cm,width = 0.5\textwidth ]{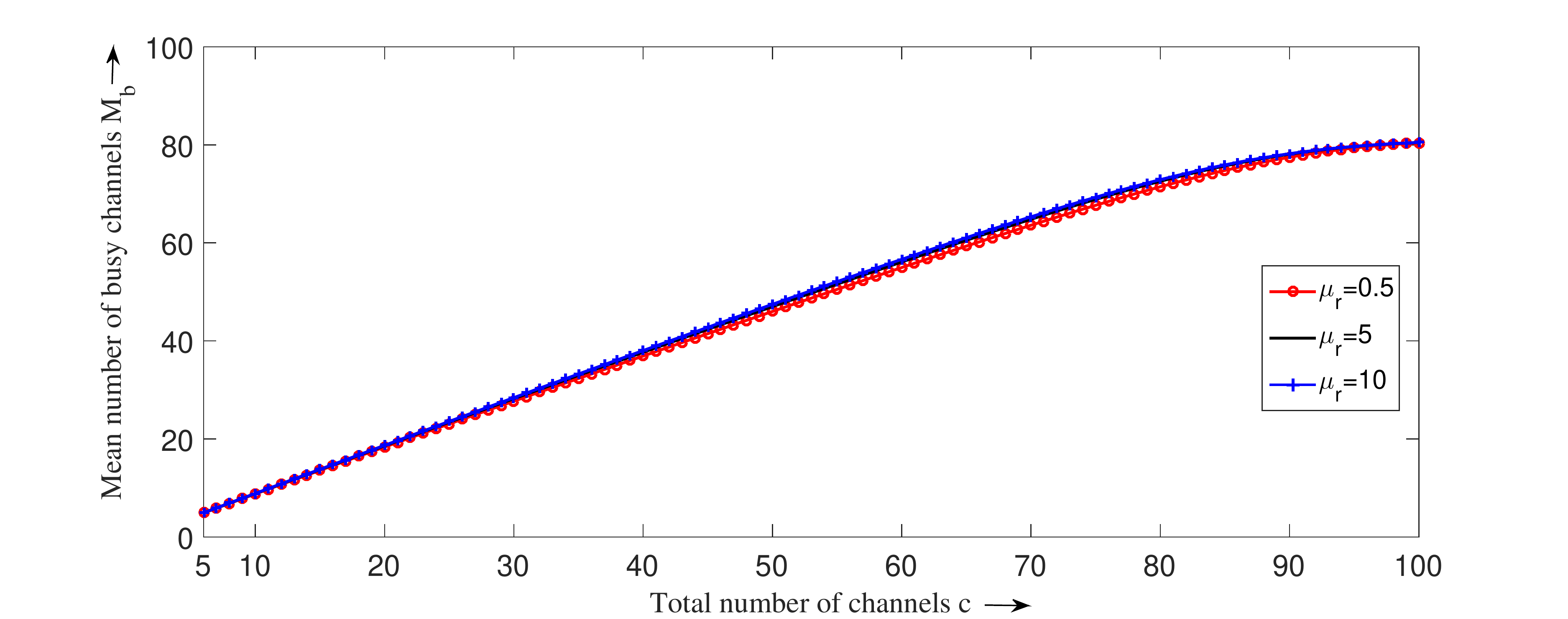}
  \label{fig:Mb_N}
}%  
\subfigure[Dependence of $M_o$ over  $c$ for $m=5$ and $g=5$.]
{\includegraphics[trim= 2cm 0.15cm 1cm 0.5cm, clip=true, height = 6cm,width = 0.52\textwidth ]{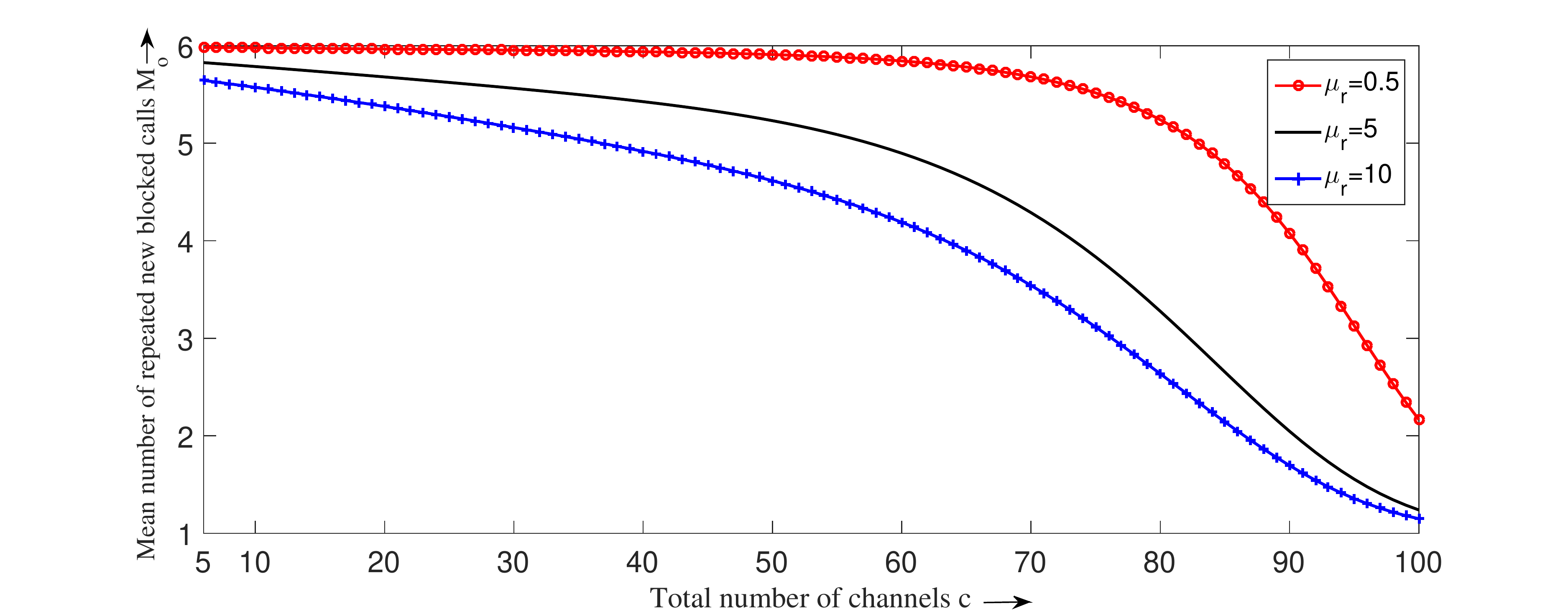} 
\label{fig:Morb_N}
}%

\label{fig:MbMorb_N}
\end{figure*}

\begin{figure}
\centering
\caption{Dependence of $M_s$ over $c$ for $m=5$ and $g=5$.}
\includegraphics[trim= 1.5cm 0.1cm 2.8cm 0.5cm, clip=true, height = 6cm,width = 0.5\textwidth]{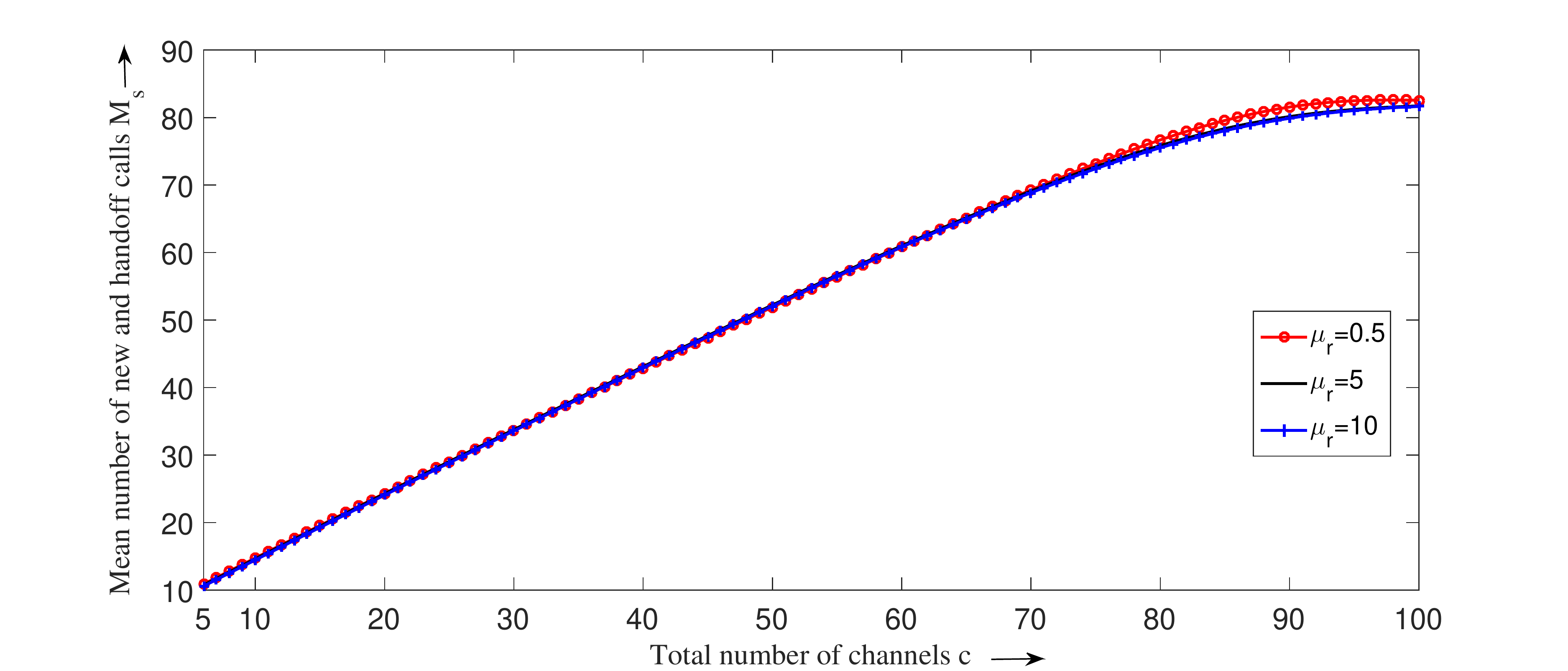}

\label{fig:Ms_N}
\end{figure}

Further, the other performance measures are discussed as follows. %Figure \ref{fig:Pbo_N} exhibits   the impact of $c$ on blocking probability of first attempt of fresh new calls $P_{bf}$. This figure presents that  $P_{bf}$ is an increasing function of $c$ for given values of   $m$ and $g$. 
Variation of  performance measures, i.e.,  mean number of busy channels $M_b$ and mean number of repeated calls in the orbit $M_{o}$, with respect to  $c$  are   shown in Figure \ref{fig:MbMorb_N}\subref{fig:Mb_N} and Figure \ref{fig:MbMorb_N}\subref{fig:Morb_N}, respectively.  It is  observed that, for fixed values of $m$ and $g$, $M_b$  increases with  $c$ and $M_{o}$   decreases with $c$. It can also be explored that   smaller values of retrial rate $\mu_{r}$ have a  negligible impact on the mean number of busy channels $M_{b}$. On the other side,    $M_o$ is a decreasing function of retrial rate $\mu_{r}$. 
Finally, Figure \ref{fig:Ms_N} exhibits the impact of $c$ on mean number of calls $M_s$. For fixed values of $m$ and $g$, $M_s$ is an increasing function of $c$. It can also be observed that for smaller values of $\mu_r$, there is a negligible change in the values of $M_s$.

%\begin{figure*}
%\centering
%\caption{Dependence of $P_{bf}$ over  $c$ for for  $g=5$.}
%{\includegraphics[trim= 1.5cm 0.1cm 2.8cm 0.6cm, clip=true, height = 6cm,width = 0.5\textwidth]{Pbo_N.pdf}

%\label{fig:Pbo_N}
%}%
%\end{figure*}

\section{Optimization Problem}
\noindent  In the proposed model,   loss probabilities are  very important performance measures.
From the customer's perspective, these performance measures need to be minimized. In this section,  by taking such factor into account, a few optimization problems are proposed to optimize total channels, guard channels and orbit size.  Specifically,    either blocking probability $P_b$ or dropping probability $P_d$ may be  minimized and  the other one may be  constrained. Thus, few  optimization problems for a multi-server retrial queueing model are presented as follows. Note that values of parameters $\lambda_{h}$, $\lambda_{n}$, $\nu$, $p$ and $\mu_r $ are considered as provided in Section 4.
\subsection*{($O1_r$)} 
Given $c$, $\lambda_{h}$, $\lambda_{n}$, $\nu$, $\mu_{r}$, $p$; calculate  the optimal integer value of $m$ and $g$ such that\\
$\begin{array}{lll}
 &\textrm{min } &P_b(m,g)\\
 &\textrm{subject to},&\\
 &&P_d(m,g) \leq P_{d0}.
 \end{array}$

\noindent Here the constant $P_{d0}$ is a pre-defined value. 

This optimization problem ($O1_r$) can be optimized by two different algorithms, provided below. 
\subsection*{Algorithm I: To find optimal orbit size.}
\begin{itemize}
\item Step-1: Set the value of  $c$. 
\item Step-2: Set $g = \lceil (x\%) (c)\rceil$, where $x$ can be assumed on the  basis of system requirement.  
\item Step-3: For such $g$, find  $m$ for which   $P_d(m, g) \leq P_{d0}$. 
%\item  It is  observed in Section 4 that  for fixed values of $c$ and $g$,  $P_b$  decreases with $m$  and $P_d$ increases with $m$.
\item Step-4:  %The largest value of $m$, which satisfies the constraint  $P_d(m, g) \leq P_{d0}$, will minimize $P_b(m)$.  
The largest value of $m$, say $m^*$ will be declared as   the optimal solution of ($O1_r$). The optimal value of ($O1_r$) will be  $P_b(m^*,g)$.
\end{itemize}
 
Optimal values for the orbit size, $m^*$, for different pre-defined values of $P_{d0}$ are provided in Table 1.

\begin{table}
\centering
\label{table:1}
\caption{Result of optimization problem ($O1_r$) by applying Algorithm I for $c=100$ and  $x=5$.}

\begin{adjustbox}{max width=0.5\textwidth}
\begin{tabular} {|c|c|c|c|c|}
\hline
$P_{d0}$ & $g$ & $m^*$ & $P_{d}(m^*, g)$ & $P_{b}(m^*, g)$ \\ [0.5ex] 
\hline
 $10^{-2}$ &5 & 0 & 0.00016136 & 0.02313149  \\  
 \hline
 $10^{-3}$ &5 & 0 & 0.00016136 & 0.02313149  \\ 
 \hline
 $10^{-4}$ &5 & 69 & 0.00009624  & 0.04732208\\
  \hline
\end{tabular}
  \end{adjustbox}

\end{table}

\subsection*{Algorithm II: To find optimal number of guard channels.}
\begin{itemize}
\item Step-1: Set the value of  $c$. 
\item Step-2: Set $m= \lceil(x\%)(c)\rceil$, where $x$ can be defined based on system requirement. 
\item Step-3: For such $m$, find  $g$ for which  $P_d(m, g) \leq P_{d0}$. 

%\item We observe in Section 4 that, for fixed values of $c$ and $m$, $P_b$  increases with $g$  and $P_d$ decreases with $g$.

\item Step-4:  The smallest value of $g$, say $g^*$ will be  declared as the optimal solution of ($O1_r$). The  optimal value will be  $P_b(m, g^*)$.
\end{itemize}
Optimal values of $g^*$ for different  values of $P_{d0}$ are mentioned in Table 2.

\begin{table}
\centering
\label{table:2}
 \caption{Result of optimization problem ($O1_r$) by applying Algorithm II for  $c=100$ and $x=5$. }

\begin{adjustbox}{max width=0.55\textwidth}
\begin{tabular} {|c|c|c|c|c|}
 \hline
 $P_{d0}$ & $m$   & $g^*$ & $P_{d}(m, g^*)$ & $P_{b}(m, g^*)$ \\ [0.5ex] 
\hline
 $10^{-2}$ & 5&  0 & 0.000786833 & 0.000784093 \\  
 \hline
 $10^{-3}$ & 5&  1 & 0.000786833 & 0.000784093 \\ 
 \hline
 $10^{-4}$ & 5  & 5 & 0.0000572980 & 0.00378360 \\
 \hline
\end{tabular}
  \end{adjustbox}
 
\end{table}

It is remarkable that if  blocked new calls do not join the orbit for retrial, i.e.,  $m=0$, the optimization problem (${O1}_r$) will be transformed  to  ($O_1$) as presented  by  \cite{haring}. Moreover,  Algorithm II will be similar to the algorithm proposed by   \cite{haring} to solve their proposed optimization problem $O_1$. Therefore, in case of no retrial, we can clearly observe  that the numerical results, summarized in Table 3, match well with results  provided by  \cite{haring}.  
\begin{table}
\centering
\label{table:3}
\caption{Result of optimization problem ($O1_r$) by applying Algorithm II for  $c=100$ and $x=0$.}

\begin{adjustbox}{max width=0.4\textwidth}
 \begin{tabular} {|c|c|c|c|}
 \hline
 $P_{d0}$ & $g^*$ & $P_{d}(g^*)$ & $P_{b}(g^*)$ \\ [0.5ex] 
 \hline
 $10^{-2}$ & 0 & 0.003992 & 0.003992 \\ 
 \hline
 $10^{-3}$ & 3 & 0.000504 & 0.012528 \\
 \hline
 $10^{-4}$ & 6 & 0.000065 & 0.023195 \\
 \hline
 $10^{-5}$ & 9 & 0.000008 & 0.038967 \\
% \hline
% $10^{-6}$ & 13 & 0.00000058 & 0.069839 \\ [1ex] 
 \hline
\end{tabular}
  \end{adjustbox}

\end{table}
%=========
\subsection*{ ($O2_r$)} 
Given $c$, $\lambda_{h}$, $\lambda_{n}$, $\nu$, $\mu_{r}$, $p$; find  the optimal integer value of $m$ and $g$ such that\\
$\begin{array}{lll}
 &\textrm{min } &P_d(m,g)\\
 &\textrm{subject to},&\\
 &&P_b(m,g) \leq P_{b0}.
 \end{array}$

\noindent Here the constant $P_{b0}$ is a pre-defined value.  The optimization problem ($O2_r$) can be optimized by the  Algorithm III.

\subsection*{Algorithm III: To find optimal orbit size}
\begin{itemize}
\item Step-1: Set  the value of  $c$. 
\item Step-2: Set $g = \lceil (x\%) (c)\rceil$, where $x$ can be defined based on system requirement. 
\item Step-3: For such $g$, find  all $m$ for which  $P_b(m,g) \leq P_{b0}$.
%\item Step-4: It is  observed in Section 4 that $P_b$ decreases and $P_d$  increases with $m$  for fixed values of $c$ and $g$.

\item Step-4: The smallest value of $m$, say $m^*$, will be declared as  the optimal solution of ($O2_r$). The  optimal value will be  $P_d(m^*,g)$.\\

\end{itemize}
Table 4 summarizes the optimal values of $m$ for different values of $P_{b0} $. 

\begin{table}
\centering
\label{table:4}
\caption{Result of optimization problem ($O2_r$) by applying Algorithm III for  $c= 105$ and $x=5$.}

\begin{adjustbox}{max width=0.45\textwidth}
 \begin{tabular} {|c|c|c|c|c|}
 \hline
 $P_{b0}$ & $g$ & $m^*$ & $P_{b}(m^*, g)$ & $P_{d}(m^*, g)$ \\ [0.5ex] 
 \hline
 $10^{-2}$ & 5& 0 & 0.0082 & 0.000046 \\
 \hline
 $10^{-3}$ & 5& 12 &0.00086  & 0.000067  \\
 \hline
 $10^{-4}$ & 5& 22 & 0.000083 & 0.000070 \\
 \hline
 \end{tabular}
  \end{adjustbox}

\end{table}

\subsection*{ ($O3_r$)}
 Given $\lambda_{h}$, $\lambda_{n}$, $\nu$, $\mu_{r}$, $p$; determine the optimal integer value of $c$, $m$ and $g$ such that\\
$\begin{array}{lll}
 &\textrm{min } &c\\
 &\textrm{subject to},&\\
 &&P_d(c,m,g) \leq P_{d0},\\
 &&P_b(c,m,g) \leq P_{b0}. 
 \end{array}$

\noindent Here constants $P_{b0}$ and $P_{d0}$ are pre-defined values.

Feasible region for ($O3_r$)  will  contain all those values of $(c,m, g)$ that satisfy both the constraints $P_d(c,m,g) \leq P_{d0}$ and $P_b(c,m,g) \leq P_{b0}$ simultaneously. This optimization problem ($O3_r$) can be optimized by the   algorithm provided below. 

\subsection*{Algorithm IV: To find optimal orbit size, optimal number of channels and optimal guard channels}
\begin{itemize}

\item Step-1: Set $m=0$.
\item Step-2: $g:=0$. Determine $c_{d0}$ and $c_{b0}$, where $c_{d0}$ is the smallest value of $c$ for which \mbox{$P_d(c,m,g) \leq P_{d0}$} and $c_{b0}$ is the smallest value of $c$ for which 
$ P_b(c,m,g) \leq P_{b0}.$ To obtain minimum value of $c$, define $c_{mid}= \displaystyle \lceil  \frac{c_{d0}+c_{b0}}{2}\rceil.$  
\item Step-3: For the value $c_{mid}$, find $g_{\max}$ and $g_{\min}$, where $g_{\max}$ is the maximum value of $g$ for which  $P_b(c_{mid},0,g) \leq P_{b0}$ and $g_{\min}$ is the minimum value of $g$ for which $P_d(c_{mid},0,g) \leq P_{d0}.$
\item Step-4: If $g_{\max}$ = $g_{\min}=g^*$, then $c_{mid}$ and $g^*$ will be optimal value of  total number of channels and optimal number of guard channels for $m=0$.
\item  Step-5: If $g_{\max} \neq g_{\min}$,  define $c_{mid}=c_{mid}+1$.  
\item Step-6: Proceed with Step 2 until  both constraints $P_d(c_{mid},0,g) \leq P_{d0}$ and $P_b(c_{mid},0,g) \leq P_{b0}$ are satisfied at $c_{mid}$ and $g^*$ . 
\end{itemize}
This procedure can  be further proceeded for $m>0$ from the Step 2.

Table 5 exhibits optimal orbit size   $m^*$, optimal number of channels  $c^*$ and optimal guard channels $g^*$ for different combinations of $P_{d0}$ and $P_{b0}$ by  applying Algorithm-IV.  We emphasize that, for $m=0$, ($O3_r$)  is converted the optimization problem $O_2$  provided by \cite{haring} and Table 5 exhibits exactly same results presented by \cite{haring}.

\begin{table}
\centering
\label{table:5}
\caption{Result of optimization problem ($O3_r$) by applying Algorithm IV.}

 \begin{adjustbox}{max width=0.65\textwidth}
 \begin{tabular} {|c|c|c|c|c|c|c|}
 \hline
$P_{d0}$ & $P_{b0}$ & $c^*$ & $m^*$ & $g^*$& $P_{b}(c^*, m^*, g^*)$ & $P_{d}(c^*, m^*, g^*)$ \\ [0.5ex]  
 \hline
 $10^{-2}$ &$10^{-1}$  & 87&0&3 &0.09089  &0.00127  \\
  \hline
 $10^{-3}$ &  $10^{-2}$&101 & 0& 2& 0.0077859 & 0.000791455 \\
  \hline
  $10^{-4}$&$10^{-3}$  &109  &0& 2&0.0009482  &0.00008555  \\
  \hline
  $10^{-5}$&$10^{-4}$ &116  &0& 2&0.0000933  &0.000007625  \\ 
  \hline
  $10^{-6}$&$10^{-5}$ &  122&0&2 &0.0000091  &0.000000687   \\[1ex] 
 \hline
\end{tabular}
 \end{adjustbox}
\end{table}

\subsection*{ ($O4_r$)}
 Given $\lambda_{h}$, $\lambda_{n}$, $\nu$, $\mu_{r}$, $p$; determine the optimal integer value of $c$, $m$ and $g$ such that\\
$\begin{array}{lll}
 &\textrm{min } &m\\
 &\textrm{subject to},&\\
 &&P_d(c,m,g) \leq P_{d0},\\
 &&P_b(c,m,g) \leq P_{b0}. 
 \end{array}$

\noindent Here constants $P_{b0}$ and $P_{d0}$ are pre-defined values. The optimization problem ($O4_r$) can be optimized by the following algorithm.

\subsection*{Algorithm V: To find optimal orbit size for given number of total channels and guard channels. }
\begin{itemize}
\item Step-1: Start with $c=2$.
\item Step-2: Set $g = \lceil (x\%) (c)\rceil$, where $x$ can be defined based on system requirement.  
\item Step-3: Find feasible region for $m$ for $m=0$ to $c/2$ such that 
  $P_d(c,m,g) \leq P_{d0}$  and  $P_b(c,m,g) \leq P_{b0}$. 
\item Step-4: Minimum value of $m$ is declared as the optimal orbit size $m^*$. Also, declare corresponding  $c$ and $g$  as  the optimal number of total channels and the optimal number of guard channels respectively. 
\item Step-5: This procedure can  be further proceeded for $c>2$.
\end{itemize}

After solving ($O4_r$) by  Algorithm V, we obtain the optimal orbit size $m^*$ corresponding to fix $c$ and $g$ for different combinations of $P_{d0}$ and $P_{b0}$ given in the Table 6.

\begin{table}
\centering
\label{table:6}
\caption{Result of optimization problem ($O4_r$) by applying Algorithm V for $x=5$.}

 \begin{adjustbox}{max width=0.45\textwidth}
 \begin{tabular} {|c|c|c|c|c|c|c|}
 \hline
$P_{d0}$ & $P_{b0}$ & $c$ & $g$& $m^*$ & $P_{b}(c,m^*,g)$ & $P_{d}(c, m^*, g)$ \\ [0.5ex]  
 \hline
 $10^{-2}$ &$10^{-1}$  &90  &4&3 &0.09089  &0.00127  \\
  \hline
 $10^{-3}$ &  $10^{-2}$& 103 &5&3 & 0.009173 & 0.000089 \\
  \hline
  $10^{-4}$&$10^{-3}$   & 112& 5 &3 &0.0007501  &0.0000061  \\
  \hline
  $10^{-5}$&$10^{-4}$ &118  &5 &3 &0.0000949  &0.00000071  \\ 
 \hline
\end{tabular}
 \end{adjustbox}

\end{table}

In practice, the total available bandwidth, distributed for new calls and handoff calls,   is limited.  By introducing the retrial phenomenon, this limited bandwidth will be shared among $c$, $g$ and $m$. Consequently, increment of orbit size $m$ will decrease the total number of channels and guard channels. This scenario is  summarized by Table 7.

\begin{table}
\centering
\label{table:7}
\caption{Impact of retrial on bandwidth sharing by $c$ and $g$.}

 \begin{adjustbox}{max width=0.5\textwidth}
 \begin{tabular} {|c|c|c|c|c|c|c|}
 \hline
$m$ &$P_{d0}$ & $P_{b0}$ & $g$ & $c$ & $P_{b}$ & $P_{d}$ \\ [0.5ex]  
 \hline
0 & $10^{-1}$ &$10^{-2}$  & 3 & 87 & 0.096834 & 0.00544\\
\hline
1 & $10^{-1}$ &$10^{-2}$  & 3 & 87 & 0.092030 & 0.005744\\
\hline
10 & $10^{-1}$ &$10^{-2}$  & 2 & 83 & 0.027259 & 0.089082\\
\hline
 \end{tabular}
 \end{adjustbox}
 
\end{table}

\section{Conclusions and Future remarks}
The main focus of this study is to   extend the work presented   by  \cite{haring}.  We then formulate  optimization problems to compute the optimal  number of channels, optimal number of guard channels and the optimal orbit size.  More specifically,  it is reported  that, for no retrial, these  optimization problems and their optimal solutions are exactly similar to the results obtained by  \cite{haring}.   Further, we    investigate numerically the impact of retrial  over the system  performance.  The blocking probability of new calls,   a  very important concept, has still not been well addressed.   Irrespective of the type of the network, i.e., 3G, 4G or 5G,  this work addresses the problem of reducing the blocking probability of new calls in cellular networks for a multi-server model with retrial phenomenon.

  We next discuss some possible extensions of this proposed model.  A closed form expression for the proposed  retrial model could be  explored in future.   Additionally,  the  direct method is adapted here to  solve $Q$-matrix and  to obtain steady state transition probabilities. This methodology can be compared with the   other methods and algorithms, e.g., matrix geometric method. Another possible extension is to consider different type of customers' arrival and service pattern.   Furthermore, it is important to consider situations that could bring the system into the inactive state. \cite{dharma2} could be a useful direction in order to study the system with such catastrophe.
  
\section*{Acknowledgements}
Authors are thankful to the editor and two anonymous reviewers for their valuable suggestions and  constructive comments  on an earlier draft. One of the authors, S. Dharmaraja,  gratefully acknowledges  the financial support received from the Department of Telecommunications (DoT), India.

\end{document}